\documentclass[leqno,twoside]{article}
\usepackage{amssymb, amsmath, latexsym, amsthm, enumerate, amsfonts}
\usepackage[sort,comma,numbers]{natbib}

\usepackage[mathscr]{euscript}
\usepackage{amsthm}
\usepackage{color}
\usepackage{amsmath}
\usepackage{mathrsfs}
\usepackage{textcomp}
\usepackage{upgreek}
\usepackage{mathabx}
\usepackage{amscd}
\usepackage{color}
\usepackage{cancel}
\usepackage{calc}

\usepackage[pdftex, pdfstartview=FitH]{hyperref}

\numberwithin{equation}{section}

\newtheorem{thm}{Theorem}[section]

\theoremstyle{definition}
\newtheorem{defn}[thm]{Definition}

\theoremstyle{remark}

\usepackage[paperwidth=165mm, paperheight=235mm, twoside, hmargin={25mm,20mm}, vmargin={20mm,15mm} ]{geometry}

\newcommand{\scal}[2]{\langle #1,#2\rangle}
\newcommand{\cc}[1]{\mathbf C^{#1}}
\newcommand{\nn}[1]{\mathbf N^{#1}}
\newcommand{\rr}[1]{\mathbf R^{#1}}
\newcommand{\mascS}{\mathscr S}
\newcommand{\maclA}{\mathcal A}
\newcommand{\maclS}{\mathcal S}
\newcommand{\maclH}{\mathcal H}
\newcommand{\mascP}{\mathscr P}
\newcommand{\op}{\operatorname{Op}}
\newcommand{\sets}[2]{\{ \, #1\, ;\, #2\, \} }
\newcommand{\Sets}[2]{\left \{ \, #1\, ;\, #2\, \right \} }

\newcommand{\cdo}{\, \cdot \, }
\newcommand{\vrum}{\vspace{0.2cm}}
\newcommand{\Sh}{\operatorname{Sh}}
\newcommand{\repart}{{\operatorname{Re}}}

\newcommand{\aw}{\operatorname{aw}}

\newcommand{\eabs}[1]{\langle #1\rangle}     
\newcommand{\nm}[2]{\Vert #1\Vert _{#2}}

\newcommand{\leqs}{\leqslant}

\begin{document}
\thispagestyle{empty}


\begin{center}
{\large \bf SOME REMARKS ON ANALYTIC PSEUDODIFFERENTIAL OPERATORS\footnote{The author would like to acknowledge the contribution of J. Toft and P. Wahlberg to the original results of the paper.}
} \vspace*{1cm}

{\bf Nenad Teofanov
\footnote{Department of Mathematics and Informatics, Faculty of Sciences, University of Novi Sad, e-mail: \href{mailto:nenad.teofanov@dmi.uns.ac.rs}{nenad.teofanov@dmi.uns.ac.rs}}
}
\end{center}


\begin{abstract}
We report some recent results on analytic pseudodifferential operators, also known as Wick operators.
An important tool in our study is the Bargmann transform which provides a coupling between the classical (real) and analytic pseudodifferential calculus. Since the Bargmann transform of Hermite functions gives rise to
formal power series in the complex domain, the results are formulated in terms of the Bargmann images of Pilipovi\'c spaces.
\\[2mm] {\it AMS Mathematics  Subject Classification $(2010)$}: 35S05, 46F05
\\[1mm] {\it Key words and phrases:}
Bargmann transform, Hermite expansions, Gelfand-Shilov spaces, Pilipovi\'c spaces, Wick and anti-Wick operators,
sharp G{\aa}rding inequality

\end{abstract}


\section{Introduction}\label{sec:intro}

We present a sample of recent results from \cite{TT2020, TTW2021} related to Wick and anti-Wick operators introduced in 1960s by Berezin in the framework of the second quantization. Wick and anti-Wick symbols are used in \cite{Berezin71} to derive various spectral properties of the corresponding operators. As demonstrated in \cite{TT2020, TTW2021}, results on Wick and anti-Wick operators provide new insight into the classical theory of pseudodifferential operators. This is done by using the mapping properties of the Bargmann transform given in \cite{Toft2012, Toft2017}.

The  Bargmann transform coupling between Hermite series expansions and formal power series expansions plays an important role in our analysis.
For that reason, we first review some facts on Hermite functions and spaces of test functions (Pilipovi\'c spaces)
with Hermite coefficients of (super)exponential decay. Thereafter we briefly discuss the  Bargmann transform, and finally we review some continuity properties of analytic  pseudodifferential operators, and sharp G{\aa}rding inequality in the context of
Wick and anti-Wick calculus.

Apart from motivations coming from quantum physics, Hermite polynomials are used e.g.
in studying the  propagation of light in infinitely long optical fibers with a parabolic index profile \cite{GT},
in visual perception and neurobiology \cite{YLM}, and in equatorial oceanography \cite{Boyd}. For the applications of pseudodifferential operators in    mobile wireless communication systems we refer to \cite{Stro}.


\section{Hermite functions and Pilipovi\'c spaces}\label{sec:hermite}

We first consider  Hermite functions  within a historical context,
and proceed with Pilipovi\'c spaces given by Hermite series with rapidly decaying coefficients.
Recall, Hermite functions are defined by
\begin{multline*}
h_\alpha (x) = \pi ^{-\frac{d}{4}}(-1)^{|\alpha |}
(2^{|\alpha |}\alpha !)^{-\frac{1}{2}}e^{\frac{|x|^2}{2}}
(\partial ^\alpha e^{-|x|^2}) \\
= e^{\frac{|x|^2}{2}} H_\alpha (x),
\;\;\;  x \in \rr d, \;\;\; n = 0,1,\dots,
\end{multline*}
where $\alpha = (\alpha_1, \alpha_2,\dots,\alpha_n) \in \nn d$,  $\alpha! = \alpha_1!\dots \alpha_d!$,
and $H_\alpha$ are (normalized) Hermite polynomials. The functions $H_n(x),$ $ x \in \mathbb{R}$, $n\in \mathbb{N}$ were
introduced by P.\ S.\ Laplace in 1810, and later studied by P.\ L.\ Chebyshev (1859) and C.\  Hermite (1864).
N. Wiener used Hermite function expansions to prove the Plancherel formula for the Fourier transform around 1930, \cite{W1933}.
In fact, Hermite functions are an orthonormal basis of $L^2 (\mathbb{R})$. Since Hermite functions are the wave functions for the stationary states of the quantum harmonic oscillator, they are  particularly useful in quantum mechanics,  \cite{F1992}.

B. Simon used Hermite function expansions in the framework of the space of rapidly decreasing functions $ \mathcal{S}(\mathbb{R}^d)$ and its dual space of tempered distributions $ \mathcal{S}'(\mathbb{R}^d)$, \cite{Simon}.
Thereafter S. Pilipovi\'c in \cite{P1988}  gave a characterization of Gelfand-Shilov type spaces and their dual spaces of tempered ultradistributions through the growth estimates of  coefficients in Hermite expansions, see also
\cite{L06,LCP} for more recent contributions in that direction.
These investigations led to a detailed study of the so-called Pilipovi\'c spaces, \cite{Toft2017, Ahmed}.

Hermite series expansions are used in the following generalization  of isotropic Gelfand-Shilov spaces
$ \mathcal{S} _{s} ( \mathbb{R}^d)$, $ s \geq 1/2$, (of Roumieu type)  which consists
of smooth functions $f$ satisfying
$$
| \partial ^\beta f (x)| \lesssim h^{|\beta|} \beta!^s e^{-k|x|^{1/s}}, \qquad x \in  \mathbb{R}^d,
$$
for some $h,k>0$. We refer to \cite{GS, Gramchev, NR, Teof2}  for details on  Gelfand-Shilov spaces and their applications in partial
differential equations. As usual, $ \mathcal{S} _{s} ' ( \mathbb{R}^d)$ denotes the dual space of the Gelfand-Shilov space
$ \mathcal{S} _{s} ( \mathbb{R}^d)$, $ s \geq 1/2$.

Pilipovi\'c spaces (of Roumieu type)  $\mathcal{H} _s (\mathbb{R}^d)$, $ s\geq 0$, are given through the formal Hermite series expansions
\begin{equation} \label{eq:hermiteexpansion}
f =\sum _{\alpha \in \mathbb{N}^d}c_\alpha h_\alpha ,\quad  c_\alpha = (f, h_\alpha), \quad |c_\alpha| \lesssim e^{-r|\alpha|^{\frac{1}{2s}}},
\end{equation}
for some $r>0$. When $ s\geq 1/2 $ we have  $\mathcal{H} _s(\mathbb{R}^ d) = \mathcal{S}  _s (\mathbb{R}^ d)$, and $\mathcal{H} _0(\mathbb{R}^d)$ is the set of all finite Hermite series expansions. It can be proved that  $\mathcal{H} _s(\mathbb{R}^d) \neq \mathcal{S} _s (\mathbb{R}^d) = \{ 0\} $, $ 1/2 > s \geq 0.$

In \cite{Toft2017} J. Toft proved that
$$
\mathcal{H} _s(\mathbb{R}^d) = \{ f \;\; | \;\; \|  \mathcal{R} ^N f \|_{L^\infty} \lesssim h^N N!^{2s} \;\; \text{for some} \;\; h >0\},
$$
where $ \displaystyle \mathcal{R} = -\frac{d^2}{ dx^2} + x^2. $ Since the Hermite functions are eigenfunctions  of $\mathcal{R} $, i.e. $\displaystyle \mathcal{R} h_n = (2n+1) h_n $, it is called the Hermite operator.

In addition, Toft considered Pilipovi\'c {\em flat} spaces where the growth condition in \eqref{eq:hermiteexpansion}
is replaced by
$$
 |c_\alpha| \lesssim r^{|\alpha| } \alpha! ^{- \frac{1}{2\sigma}}, \;\;\; \sigma > 0
$$
some $r>0$.

The well known relation between $ L^2(\mathbb{R}^d) $
and the Fock space of analytic functions $A^2(\cc d) $ (see Section \ref{sec:bargmann} for the definition)
can then be extended to the relation between Pilipovi\'c spaces and specific subspaces of the space of analytic functions,
\cite{Toft2012}. This is done via the Bargmann transform, cf. Definition \ref{def:bargmann}.
Following this approach, a detailed study of analytic pseudodifferential operators
is given in \cite{TT2020, TTW2021}.


\section{The Bargmann transform}\label{sec:bargmann}

\begin{defn} \label{def:bargmann}
The Bargmann transform $\mathfrak V_df$ of $f\in \maclS _{1/2}'(\rr d)$
is the entire function
\begin{multline*}
\mathfrak V_d f (z) =\int_{\rr d} \mathfrak A_d(z,y)f(y)\, dy
\\[1ex]
=
\pi ^{-\frac d4}\int _{\rr d}\exp \Big ( -\frac 12(\scal
z z+|y|^2)+2^{1/2}\scal zy \Big )f(y)\, dy,
\end{multline*}
$ z \in \cc d$, and $ \scal zw = \sum _{j=1}^dz_jw_j $.
\end{defn}

It was proved by V. Bargmann in 1961. that
$$
\mathfrak V_d : L^2(\rr d) \to A^2(\cc d)
$$
is a bijective and isometric mapping from $L^2(\rr d)$ to the Fock space $A^2(\cc d) $,
the Hilbert space of entire functions with the scalar product
\begin{equation*}
(F,G)_{A^2}\equiv  \int _{\cc d} F(z)\overline {G(z)}\, d\mu (z),\quad F,G\in A^2(\cc d),
\end{equation*}
where $d\mu (z)=\pi ^{-d} e^{-|z|^2}\, d\lambda (z)$ ($d\lambda (z)$ is
the Lebesgue measure on $\cc d$).

These investigations put a solid theoretical background for a quantization procedure proposed by V. Fock back in 1930's.
More precisely,  Bargmann showed that $\mathfrak V_d $
maps the creation and annihilation operators,
$ A = -\frac{d}{ dx} + x, $ and $ A^{\dag} = \frac{d}{ dx} + x$ respectively,
into multiplication and differentiation in the complex domain,
\cite{B1, B2}. Note that $ \mathcal{R} = \frac{1}{2} (A A^{\dag} + A^{\dag}  A)$.

The Bargmann transform maps the Hermite functions to monomials as
\begin{equation*}\label{Eq:BargmannHermiteMap}
\mathfrak V_d h_\alpha = e_\alpha ,\qquad e_\alpha (z)= \frac {z^\alpha}{\alpha !^{\frac 12}},
\quad z\in \cc d,\quad \alpha \in \nn d.
\end{equation*}

The orthonormal basis
$\{ h_\alpha \}_{\alpha \in \nn d} \subseteq L^2(\rr d)$
is thus mapped to the orthonormal basis
$\{ e_\alpha \} _{\alpha \in \nn d}\subseteq A^2(\cc d)$.

Let  $\maclA _0(\cc d) $ be the set of all analytic polynomials of the form $F (z) =  \sum_{|\alpha| \leq N } c(F,\alpha) e_\alpha (z)$, for some $N\in \mathbb{N}$, and let
$$
\maclA _s(\cc d) = \{ F (z) = \sum_{\alpha \in \mathbb{N} ^d} c(F,\alpha) e_\alpha (z) \;\; | \;\; | c(F,\alpha)| \lesssim e^{-r|\alpha|^{\frac{1}{2s}}} \}, \;\;\; s>0.
$$
Then it is proved by Toft (\cite{Toft2017}) that
$$ \mathfrak V_d :  \maclH _s(\rr d) \rightarrow \maclA _s(\cc d), \qquad s>0, $$
is bijective mapping between Pilipovi\'c spaces and corresponding spaces of analytic functions.


\section{Analytic pseudodifferential operators}\label{sec:psido}

\begin{defn}
Let $a $ be a locally bounded function on $\cc {2d}$ such that
$(z,w) \mapsto a (z,w)$ is analytic, $z,w\in \cc d$. The {\em analytic pseudodifferential operator} or the Wick operator
$\op _{\mathfrak V}(a)$  with {\em the symbol} $a$ is given by
$$
\op _{\mathfrak V}(a)F (z) =
\pi ^{-d}
\int _{\cc d} a (z,w)F(w)e^{(z-w,w)}\, d\lambda (w),
$$
where $F$ is an entire function, $d\lambda$ is the Lebesgue measure
and $(\cdo ,\cdo )$ is the scalar product on $\cc d$.
\end{defn}

Thus  $(\op _{\mathfrak V}(a)F)(z)$ is equal to the integral
operator
$$
(T_KF)(z) =
\pi ^{-d} \int _{\cc d} K(z,w)F(w)\,  e^{-|w|^2}\, d\lambda (w)
=
\int _{\cc d} K(z,w)F(w)\, d\mu (w),
$$
when $K(z,w)=$ $ K_a(z,w) = a(z,w)e^{(z,w)},$
and  $d\mu (w)=\pi ^{-d} e^{-|w|^2}\, d\lambda (w)$.

The (classical) pseudodifferential operator $\op (b) $ with the symbol  $b$
is given by the Kohn-Nirneberg correspondence
$$
f(x)\mapsto (\op (b)f)(x) = (2\pi )^{-\frac d2}\int _{\rr d}
b(x,\xi )\widehat f(\xi )e^{i\scal x\xi}\, d\xi.
$$

If $b$ is a polynomial symbol, i.e.
$$ b(x,\xi) = \sum_{|\alpha +\beta| \leq N} c_1 (\alpha, \beta) x^\alpha \xi^\beta,$$
then there is a unique symbol
$$ a(z,w ) = \sum_{|\alpha +\beta| \leq N} c_2 (\alpha, \beta) z^\alpha \overline{w}^\beta$$
such that
$ \displaystyle \op _{\mathfrak V}(a) = \mathfrak V_d \circ \op (b) \circ \mathfrak V_d ^{-1}$
(cf. \cite{TTW2021}).

Let $ \wideparen \maclA _{s}(\cc {2d})
= \sets {K}{(z,w)\mapsto K(z,\overline w)\in \maclA _{s}(\cc {2d})} $, $ s \geq 0,$
and $$
\wideparen A(\cc {d_1 + d_2})
\equiv
\Sets {K(z,\overline w),\quad z\in \cc {d_1}, w\in  \cc {d_2}}
{K \;\; \text{ is an analytic function} \; }.
$$

Identification of linear and continuous mappings with pseudodifferential operators, and their
basic continuity properties are given in Theorems \ref{Thm:GSPseudoDifficultDirection} and \ref{Thm:GSPseudoEasyDirection}.
We refer to \cite{TT2020} for the proofs.

\begin{thm}\label{Thm:GSPseudoDifficultDirection}
Let $s\ge \frac 12$. Then the following is true:
\begin{enumerate}
\item If $T$ is a linear and continuous map
from $\maclA _s'(\cc {d})$ to
$\maclA _s(\cc {d})$, then there is a unique
$a\in \wideparen A (\cc{d}\times \cc {d})$ such that
$$
|a(z,w)|\lesssim e^{\frac 12\cdot |z-w|^2-r(|z|^{\frac 1s}+|w|^{\frac 1s})},
\quad z,w \in \cc {d},
$$
for some $r>0$ and $T =\op _{\mathfrak V}(a)$;

\vrum
\item If $T$ is a linear and continuous map
from $\maclA _s(\cc {d})$ to
$\maclA _s'(\cc {d})$, then there is a unique
$a\in \wideparen A (\cc{d}\times \cc {d})$ such that
$$
|a(z,w)|\lesssim e^{\frac 12\cdot |z-w|^2+r(|z|^{\frac 1s}+|w|^{\frac 1s})},
\quad z,w \in \cc {d},
$$
for every  $r>0$ and $T =\op _{\mathfrak V}(a)$.
\end{enumerate}
\end{thm}

\begin{thm}\label{Thm:GSPseudoEasyDirection}
Let $s\ge \frac 12$. Then the following is true:
\begin{enumerate}
\item If $a\in \wideparen A(\cc {d}\times \cc {d})$ satisfies
$$
|a(z,w)|\lesssim e^{\frac 12\cdot |z-w|^2-r(|z|^{\frac 1s}+|w|^{\frac 1s})},
\quad z,w \in \cc {d},
$$
for some $r>0$,
then $\op _{\mathfrak V}(a)$ from $\maclA _0(\cc d)$ to $\maclA _0'(\cc d)$
is uniquely extendable to a linear and continuous map from
$\maclA _s'(\cc {d})$ to $\maclA _s(\cc {d})$;

\vrum

\item If $a\in \wideparen A(\cc {d}\times \cc {d})$ satisfies
$$
|a(z,w)|\lesssim e^{\frac 12\cdot |z-w|^2+r(|z|^{\frac 1s}+|w|^{\frac 1s})},
\quad z,w \in \cc {d},
$$
for every  $r>0$,
then $\op _{\mathfrak V}(a)$ from $\maclA _0(\cc d)$ to $\maclA _0'(\cc d)$
is uniquely extendable to a linear and continuous map from
$\maclA _s(\cc {d})$ to $\maclA _s'(\cc {d})$ .
\end{enumerate}
\end{thm}

An important subclass of Wick operators are the
anti-Wick operators, which are
Wick operators where the symbol
$a(z,w)$ does not depend on $z$:
$$
\op _{\mathfrak V}^{\aw}(a_0)F (z) = \pi ^{-d}
\int _{\cc d} a_0(w)F(w)e^{(z-w,w)}\, d\lambda (w).
$$
The anti-Wick
operators can also be described as the Bargmann images of
Toeplitz operators on $\rr d$.
We refer to \cite{LieSol, Sh1, Toft2012} for more details, and note that an important feature in energy
estimates in quantum mechanics and
time-frequency analysis is that  non-negative symbols give rise to non-negative Toeplitz  and anti-Wick operators.

In the next theorem we show that many Wick operators can essentially be expressed
as linear combinations of anti-Wick operators. The expansion
\eqref{Eq:WickToAntiWick} is deduced by using
Taylor expansion and integration by parts, see \cite{TTW2021} for details.

\begin{thm} \label{thm:wick-antiwick}
Suppose $s\ge \frac 12$, $a \in \wideparen \maclA _s'(\cc {2d})$ (the dual of $\maclA _s (\cc {2d})$),
let $N\ge 1$ be an integer, and let
\begin{align*}
a_\alpha (w) &= \partial _z^\alpha \overline \partial _w^\alpha a (w,w), \quad \alpha \in \nn d,
\intertext{and}
b_\alpha (z,w) &= |\alpha |
\int _0^1 (1-t)^{|\alpha |-1} \partial _z^\alpha \overline \partial _w^\alpha a
(w+t(z-w),w)\, dt, \quad \alpha \in \nn d\setminus 0.
\end{align*}
Then
\begin{equation}\label{Eq:WickToAntiWick}
\op _{\mathfrak V}(a)
=
\sum _{|\alpha |\le N}
\frac {(-1)^{|\alpha |}\op _{\mathfrak V}^{\aw}(a_\alpha )}{\alpha !}
+
\sum _{|\alpha | = N+1}\frac {(-1)^{|\alpha |}\op _{\mathfrak V}(b_\alpha )}{\alpha !}.
\end{equation}
\end{thm}

We apply Theorem \ref{thm:wick-antiwick} to the sharp  G{\aa}rding inequality for analytic pseudodifferential operators.
Its real counterpart represents one of the basic applications of the Anti-Wick theory, \cite{NR}.
According to G.\ Folland, G{\aa}rding's inequality is a milestone in the theory of elliptic equations, \cite{Fol1995}.

\vspace{3mm}

As a preparation, we introduce the Shubin class of symbols (cf. \cite{Sh1}).
First we introduce weight functions as follows.
A weight on $\rr d$ is a positive function $\omega \in  L^\infty _{loc}(\rr d)$
such that $1/\omega \in  L^\infty _{loc}(\rr d)$. It
is $v-$moderate  for a polynomially bounded weight if there is another weight
$v$ of the form $v(x) = {\eabs x}^s$, $s \ge 0$
($\eabs x = (1+|x|^2)^{\frac 12}$) such that
\begin{equation*}
\omega(x+y)\lesssim \omega(x)v(y),\quad x,y\in\rr{d}.
\end{equation*}
By  $\mascP _{\Sh ,\rho }(\rr d)$, $0 \leq \rho \leq 1$, we denote
the set of all smooth and $v-$moderate weights $\omega $
for a polynomially bounded weight $v$ such that
$$
|\partial ^\alpha \omega (x)|\lesssim \omega (x)\eabs x^{-\rho |\alpha |},
\quad \alpha \in \nn d, \;\; x \in \rr d.
$$

Let $0 \leq \rho \leq 1$, and let $\omega \in \mascP _{\Sh ,\rho }(\rr {d})$. The Shubin symbol class
$\Sh _\rho ^{(\omega )}(\rr {d})$ is the set of all $a\in C^\infty (\rr {d})$
such that
$$
|\partial ^\alpha a(x)|\lesssim \omega (x) \eabs x^{-\rho |\alpha |},
\qquad x \in \rr {d},
$$
for every multi-index $\alpha \in \nn {d}$.

Let $\wideparen \maclA _{\Sh ,\rho}^{(\omega )}(\cc {2d})$,
be the set of all $a\in \wideparen A(\cc {2d})$ such that
\begin{equation}\label{Eq:ShubinWickEstimates}
\left| \partial_z^\alpha \overline \partial _w
^\beta a(z,w) \right |
\le C
e^{ \frac1{2} |z-w|^2}\omega (\sqrt 2\, \overline z)\eabs{z+w}^{- \rho |\alpha + \beta|}
\eabs{z-w}^{-N} ,\quad N\ge 0.
\end{equation}

Let $a_0 \in \maclS _{1/2}'(\rr {2d})$.
Then the Bargmann assignment $\mathsf S_{\mathfrak V} a_0$ of $a_0$ is
the unique element $a\in \wideparen A (\cc {2d})$ which fulfills
$$
\op _{\mathfrak V}(a) = \mathfrak V_d \circ \op^w (a_0 )\circ \mathfrak V_d^*
\quad \Leftrightarrow \quad a = \mathsf S_{\mathfrak V} a_0,
$$
where $\op^w (a_0 )$ is  the Weyl pseudodifferential operator
$$
\op^w (a_0 )f (x)
=
(2\pi  ) ^{-d}\iint a_0 (\frac{x+y}{2}),\xi )
f(y)e^{i\scal {x-y}\xi }\,
dyd\xi, \quad x\in \rr d.
$$
It can be proved that  ${\mathsf S}_{\mathfrak V}$ is a homeomorphism from
$\Sh ^{(\omega )}_\rho (\rr {2d})$ to $\wideparen
\maclA _{\Sh ,\rho}^{(\omega )}(\cc {2d})$,  $0 \le \rho \leqs 1$, $\omega \in \mascP _{\Sh ,\rho}(\rr {2d})$, see \cite{TTW2021}.

\vspace{3mm}

Finally, we have the following version of the sharp  G{\aa}rding inequality.

\begin{thm}\label{Thm:ShGarding}
Let $\rho >0$, $\omega (z) = \eabs z ^{2\rho}$ and let
$a_0\in \wideparen \maclA _{\Sh ,\rho}^{(\omega )}(\cc {2d})$ be such that
$a_0(w,w)\ge -C_0$ for all $w \in \cc d$, for some constant $C_0\ge 0$. Then
\begin{alignat*}{2}
\repart ((\op _{\mathfrak V}(a_0)F,F)_{A^2})
&\ge
-C\nm F{A^2}^2, &
\qquad
F &\in
{\mathfrak V}_d (\mascS (\rr d))
\intertext{and}
|\operatorname{Im}((\op _{\mathfrak V}(a_0)F,F)_{A^2})|
&\le
C\nm F{A^2}^2, &
\qquad
F &\in
{\mathfrak V}_d (\mascS (\rr d))
\end{alignat*}
for some constant $C\ge 0$.
\end{thm}

We refer to \cite{TTW2021} for the proof.


\section*{Acknowledgement}

This work is partially supported by project DS15 TIFREFUS, MPNTR Grant No. 451-03-9/2021-14/200125,
and Project 19.032/961-103/19 MNRVOID of the Republic of Srpska.


\end{document}